\documentclass[12pt]{article}
\usepackage{amssymb}
\usepackage{amsmath}

\newtheorem{theorem}{Theorem}

\newtheorem{example}{Example}

\newcommand{\gap}{\vspace{0.1in}}

\newcommand{\goes}{\rightarrow}

\newcommand{\rnn}{R^{n\times n}}

\newcommand{\qed}{\Box}

\newcommand{\Rn}{R^n}

\newcommand{\Rnplus}{R_+^n}
\newcommand{\Sn}{{\cal S}^n}
\newcommand{\Snplus}{{\cal S}^n_+}
\newcommand{\Ln}{{\cal L}^n}
\newcommand{\Lnplus}{{\cal L}^n_+}

\title{\bf
A representation theorem for the Lorentz cone automorphisms 
    }
\author{
Roman Sznajder\\
Department of Mathematics\\
Bowie State University\\
Bowie, MD 20715\\rsznajder@bowiestate.edu
}
\date{February 15, 2021}
\begin{document}
\maketitle
\begin{abstract}
\noindent In this note we prove a representation theorem for the symmetric cone automorphisms in the spin algebra\, $\Ln$. 
\end{abstract}

\vspace{1cm}




\section{Introduction}

Euclidean Jordan algebras have become a subject of intensive studies in the current literature on conic optimization (see, e.g., Schmieta and Alizadeh \cite{schmieta-alizadeh}). Euclidean Jordan algebra (EJa) is a finite dimensional real inner product space \,$V$\, with a bilinear mapping \,$(x,y) \goes x \circ y$\, satisfying certain properties that will be described in the next section. The spin algebra \,$\Ln$ $(n>1$), which is a specific example of EJa, has been heavily studied in relation to the so-called second-order cone optimization, see, e.g., Fukushima, Luo, and Tseng \cite{fukushima etal}.  

In any Euclidean Jordan algebra \,$V$\, there is the cone of squares \,$K=\{x \circ x: x \in V\}$\, which is self-dual and homogeneous cone convex cone. In such an algebra one can define the automorphism groups \,$Aut(V)$\, and \,$Aut(K)$\, the following way (see Faraut and Kor\'{a}nyi \cite{faraut-koranyi}): $S$ \, is an {\em algebra automorphism} ($S \in Aut(V)$)\, if \,$S:V \goes V$\, is an invertible linear transformation satisfying the condition \,$S(x \circ y)=S(x) \circ S(y)$\, for all \,$x, y \in V$, and \,$S$\, is a {\em cone automorphism} ($S \in Aut(K)$)\, if \,$S$\, is a linear transformation satisfying the condition\, $S(K)=K$.

There is a vast knowledge on groups of algebra and cone automorphisms in a general setting of EJas (Faraut and Kor\'{a}nyi \cite{faraut-koranyi}), but only in several instances an intrinsic characterization of such automorphisms can be given. It is known in the case of algebra \,$\Rn$,\, algebra of all real symmetric matrices \,$\Sn$ (Schneider \cite{schneider}), and \,$\Ln$\, (Loewy and Schneider \cite{loewy-schneider}). Our aim is to provide a constructive characterization of the cone automorphism group \,$Aut(\Ln_+)$. 
  

\section{Preliminaries}
In this section, we briefly recall the definition and give some examples of EJas.
The relevant material can be found in \cite{faraut-koranyi}. For connections with complementarity theory see Gowda, Song \cite{gowda-song}, Gowda, Sznajder, and Tao \cite{gowda-sznajder-tao}, as well as Gowda and Sznajder \cite{gowda-sznajder}.

A  {\it Euclidean Jordan algebra} is a triple $(V, \circ, \langle\cdot, \cdot\rangle)$ where $(V,
\langle\cdot,
\cdot\rangle)$ is a finite dimensional
inner product space over $R$ and $(x,y)\mapsto x\circ y:V\times V\to V$ is a bilinear mapping
satisfying the following conditions:

\begin{itemize}
\item [(i)] $x\circ y=y\circ x$ for all $x,y\in V$,
\item [(ii)] $x\circ(x^2\circ y)=x^2\circ (x\circ y)$ for all $x,y\in V$ where $x^2:=x\circ x$, and
\item [(iii)] $\langle x\circ y,z\rangle =\langle y, x\circ z\rangle$ for all $x,y,z\in V.$
\end{itemize}
In addition, we assume that there is an element $e\in V$ (called the {\it unit} element) such that $x\circ e=x$ for all
$x\in V$.

Here are several examples of Euclidean Jordan algebras and the corresponding automorphism groups.

\begin{example} \label{example 1} Consider $R^n$ with the usual inner product $\langle\cdot,\cdot\rangle$ and the
componentwise product defined by $(x*y)_i=x_iy_i$ for $i=1,2,\ldots,n$, where $x_i$ is the $i$th component of (column) vector $x\in R^n$. The corresponding cone of squares is $K=\Rnplus$, the nonnegative orthant in $\Rn$. One can easily see that $Aut(R^n)$ consists of permutation matrices and any element in Aut$(\Rnplus)$ is  a product of a permutation matrix and a diagonal matrix with positive diagonal entries.
\end{example}

\begin{example} \label{example 2} Let $\Sn$ be the set of all $n\times n$ real symmetric matrices with the inner and
Jordan product given by
$$\langle X,Y\rangle:=trace(XY)\quad\mbox{and}\quad X\circ Y:=\frac{1}{2}(XY+YX).$$
In this setting, the cone of squares $\Snplus$ is the set of all positive semidefinite matrices in $\Sn$. 

It is known (see Schneider \cite{schneider}) that the corresponding to any
$\Gamma \in Aut(\Sn_+)$, there exists an invertible matrix $Q \in \rnn$ such that
$$ \Gamma(X) =  QXQ^T~~ (X \in \Sn).$$
Also for $\Lambda\in Aut(\Sn)$, there exists an orthogonal matrix $U$ such that 
$$\Lambda(X)=UXU^T~~ (X \in \Sn).$$ 
\end{example}

\vspace{-.8cm}
\begin{example} \label{example 3} Consider $R^n$ $(n>1)$ where any element $x$ is written as
$$x=\left [ \begin{array}{c}
x_0\\\overline{x}\end{array}\right ]$$
with $x_0\in R$ and $\overline{x}\in R^{n-1}$. The inner product in $\Rn$ is the usual inner product.
The Jordan product $x\circ y$ in $R^n$ is defined by
$$
x\circ y=\left [ \begin{array}{c}
x_0\\
\overline{x}\end{array}\right ]\circ \left [ \begin{array}{c}
y_0\\
\overline{y}\end{array}\right ]:=\left [ \begin{array}{c}
\langle x,y\rangle\\
x_0\overline{y}+y_0\overline{x}\end{array}\right ].$$
We shall denote this Euclidean Jordan algebra
$(R^n,\circ,\langle\cdot,\cdot\rangle)$ by $\Ln$.
In this algebra, the cone of squares, denoted by $\Lnplus$, is called the {\it Lorentz cone}
(or the second-order cone). It is
given by
$$\Lnplus=\{x: ||\overline{x}||\leq x_0\}.$$ 

In this case, it is known (see Loewy and Schneider \cite{loewy-schneider}) that if $n \times n$
matrix $\Gamma$ belongs
to $Aut({\cal L}^n_+)$, then there exists $\mu >0$ such that
$$\Gamma^TJ_n\Gamma= \mu J_n\,,$$
where $J_n=diag(1, -1, -1, \ldots ,-1)$. In particular, if $\Lambda \in Aut({\cal L}^n)$, then 
(because $\Lambda(e)=e$), it can be easily
seen that
$$\Lambda= \left [\begin{array}{cc} 1 & 0\\0 & D \end{array} \right ]$$
where $D: R^{n-1} \rightarrow R^{n-1}$ is an orthogonal matrix.
\end{example}


\section{A representation theorem for the Lorentz cone\\ automorphisms}
In this section we will prove the following representation theorem
\begin{theorem} \label{repr-theorem}
Let $S$ be an $n \times n$ real matrix. Then $S \in Aut(\Ln_+)$ if and only if there exist $\nu >0$, $\alpha \geq 0$, and $U,\, V$ orthogonal transformations in $R^{n-1}$
such that
\begin{equation} \label{repr}
S= \nu \left [ \begin{array}{cc} 1 & 0\\0 & V \end{array} \right ]
\left [ \begin{array}{cc|c} \sqrt{1+\alpha^2} & \alpha & {\bf 0}^T\\ 
\alpha &  \sqrt{1+\alpha^2} & \\
\hline 
{\bf 0} & & I_{n-2} \end{array} \right ]
\left [ \begin{array}{cc} 1 & 0\\0 & V^T \end{array} \right ]
\left [ \begin{array}{cc} 1 & 0\\0 & U \end{array} \right ]
\end{equation}
and conversely.
\end{theorem}

\noindent {\bf Proof.} Let $S \in Aut(\Ln_+)$. We know that there exists $\mu >0$ such that $S^TJS=\mu J=SJS^T$, see Loewy and Schneider \cite{loewy-schneider}. By scaling $S$ appropriately, we assume that $\mu =1$. Thus, 
\begin{equation} \label{scaled Loewy-Schneider} 
S^TJS=J=SJS^T.
\end{equation}
Let $S=\left [ \begin{array}{cc} a & b^T\\c & D \end{array} \right ]$. Our first observation is that $a > ||c||$, as
$Se=\left [ \begin{array}{c} a \\ c \end{array} \right ] > 0$. Using the first equality in (\ref{scaled Loewy-Schneider}),
we have (with $I=I_{n-1}$)
$$\left [ \begin{array}{cc} a & c^T\\b & D^T \end{array} \right ]
\left [ \begin{array}{cc} 1 & 0 \\0 & -I \end{array} \right ] \left [ \begin{array}{cc} a & b^T\\c & D\end{array} \right ]=
\left [ \begin{array}{cc} 1 & 0 \\0 & -I \end{array} \right ].$$
This results in
$$\left [ \begin{array}{cc} a^2=||c||^2 & ab^T-c^TD \\ab-D^Tc & bb^T-D^TD \end{array} \right ]=
\left [ \begin{array}{cc} 1 & 0 \\0 & -I \end{array} \right ].$$

\gap
\noindent Then
$$ \begin{array}{ll} (A1)~~ & a=\sqrt{1+||c||^2} \hspace{3.5in}\\(A2) & ab=D^Tc\\(A3) & D^TD=I+bb^T .\end{array}  $$

Using the second equality in (\ref{scaled Loewy-Schneider}), we get
$$ \begin{array}{ll} (B1)~~ & a=\sqrt{1+||b||^2} \hspace{3.5in}\\ (B2) & ac=Db\\(B3) &  DD^T=I+cc^T .\end{array} $$

By the polar decomposition, $D=PU$, with $P$ is symmetric, positive definite and $U$ orthogonal. Then $D^T=U^TP^T=U^TP$,
$DD^T=PUU^TP=P^2$, so $P=\sqrt{DD^T}$. Hence, 
\begin{equation} \label{auto-repr1}
S=\left [ \begin{array}{cc} a & \frac{1}{a}(D^Tc)^T\\c & PU \end{array} \right ]=
\left [ \begin{array}{cc} a & \frac{1}{a}c^TPU\\c & PU \end{array} \right ]=
\left [ \begin{array}{cc} a & \frac{1}{a}c^TP\\c & P \end{array} \right ]
\left [ \begin{array}{cc} 1 & 0\\0 & U \end{array} \right ].
\end{equation}
Actually, assuming (A1), (A2), and (B3), the other properties come out automatically. 

\gap
\noindent To show (B2): Since $DD^T=I+cc^T$, $DD^Tc=(I+cc^T)c=c+||c||^2c=a^2c$. Thus, $aDb=DD^Tc=a^2c$, so $Db=ac$.

\gap
\noindent To show (A3): By using (B2), we have 
$$(I+bb^T)D^T=D^T+b(Db)^T=D^T+b(ac^T)=D^T+abc^T.$$
By using (B3) and (A2), we get
$$D^TDD^T=D^T(I+cc^T)=D^T+D^Tcc^T=D^T+abc^T.$$
Altogether, $(I+bb^T)D^T=D^TDD^T$, so $I+bb^T=D^TD$.

\gap
\noindent To show (B1): By using (B3) and (A3)
$$(n-1)+||c||^2=tr(DD^T)=tr(D^TD)=(n-1)+||b||^2.$$
Hence, $||b||=||c||$, so (B1). 

\gap
For $c \in R^{n-1}$, let $V$ be an orthogonal matrix and a real number $\alpha$ such that
$c=\alpha V e_1$, where $e_1=[1,0, \cdots, 0]^T \in R^{n-1}$. By (B3),
$$P=\sqrt{DD^T}=\sqrt{I+cc^T}=\sqrt{I+\alpha^2Ve_1e_1^TV^T}=V\sqrt{I+\alpha^2e_1e_1^T} V^T.$$
After substituting in (\ref{auto-repr1}),
$$S=\left [ \begin{array}{cc} \sqrt{1+\alpha^2} & \frac{1}{\sqrt{1+\alpha^2}}\alpha e_1^TV^TV\sqrt{I+\alpha^2e_1e_1^T} V^T\\ \mbox{} \vspace{-.3cm}&  \\ 
\alpha Ve_1 &  V\sqrt{I+\alpha^2e_1e_1^T} V^T\end{array} \right ]
\left [ \begin{array}{cc} 1 & 0\\0 & U \end{array} \right ]=$$
$$\left [ \begin{array}{cc} 1 & 0\\0 & V \end{array} \right ]
\left [ \begin{array}{cc} \sqrt{1+\alpha^2} & \frac{1}{\sqrt{1+\alpha^2}}\alpha e_1^T\sqrt{I+\alpha^2e_1e_1^T}\\
 \mbox{} \vspace{-.3cm}&  \\  
\alpha e_1 &  \sqrt{I+\alpha^2e_1e_1^T} \end{array} \right ]
\left [ \begin{array}{cc} 1 & 0\\0 & V^T \end{array} \right ]
\left [ \begin{array}{cc} 1 & 0\\0 & U \end{array} \right ].$$

\gap
\noindent \underline{Claim}: $\sqrt{I+\alpha^2e_1e_1^T} e_1=ae_1$.\\
In fact, let $E:=\sqrt{I+\alpha^2e_1e_1^T}$, so $E^2e_1=(I+\alpha^2e_1e_1^T)e_1=e_1+\alpha^2e_1=a^2e_1$. Thus
$$(E+aI)(E-aI)e_1=0.$$
Since $E+aI$ is invertible (as a positive definite matrix), $(E-aI)e_1=0$. \\
After the simplification given in Claim,
$$S=\left [ \begin{array}{cc} 1 & 0\\0 & V \end{array} \right ]
\left [ \begin{array}{cc} \sqrt{1+\alpha^2} & \alpha e_1^T\\ 
\alpha e_1 &  \sqrt{I+\alpha^2e_1e_1^T} \end{array} \right ]
\left [ \begin{array}{cc} 1 & 0\\0 & V^T \end{array} \right ]
\left [ \begin{array}{cc} 1 & 0\\0 & U \end{array} \right ].$$

\gap
\noindent Since 
$$e_1e_1^T=\left [\begin{array}{cccc} 1 & 0 & \cdots & 0\\
0 & & & 0\\
\vdots & & & \vdots\\
0 & 0 & \cdots & 0 \end{array} \right ],$$
we finally get
$$S=\left [ \begin{array}{cc} 1 & 0\\0 & V \end{array} \right ]
\left [ \begin{array}{cc|c} \sqrt{1+\alpha^2} & \alpha & {\bf 0}^T\\ 
\alpha &  \sqrt{1+\alpha^2} & \\
\hline 
{\bf 0} & & I_{n-2} \end{array} \right ]
\left [ \begin{array}{cc} 1 & 0\\0 & V^T \end{array} \right ]
\left [ \begin{array}{cc} 1 & 0\\0 & U \end{array} \right ].$$

\gap
Conversely, a direct computation shows that the transformation 
$$T_{\alpha}:=\left [ \begin{array}{cc|c} \sqrt{1+\alpha^2} & \alpha & {\bf 0}^T\\ 
\alpha &  \sqrt{1+\alpha^2} & \\
\hline 
{\bf 0} & & I_{n-2} \end{array} \right ]$$
is a cone automorphism of $\Ln_+$, so any linear transformation of the form (\ref{repr}) 
is a cone automorphism of $\Ln_+$ as well.
\mbox{} \hfill $\qed$

\noindent {\bf Remark.} Note that any $S \in Aut(\Ln_+)$ can be represented as
\begin{equation} \label{auto-repr2}
S=\nu \left [ \begin{array}{cc} a & c^T\\c & P \end{array} \right ]
\left [ \begin{array}{cc} 1 & 0\\0 & U \end{array} \right ],
\end{equation}
where $\nu >0$, $c \in R^{n-1}$, $a=\sqrt{1+||c||^2}$,  $P=\sqrt{I+cc^T}$, and $U$---orthogonal transformation in $R^{n-1}$.
In fact, according to (\ref{auto-repr1}) (when $\nu=1)$, 
$$S=\left [ \begin{array}{cc} a & \frac{1}{a}c^TP\\c & P \end{array} \right ]
\left [ \begin{array}{cc} 1 & 0\\0 & U \end{array} \right ].$$
Since $P^2=I+cc^T$, $P^2c=(1+||c||^2)c=a^2c$. Thus, $Pc=ac$ and (\ref{auto-repr2}) has been proved. \hfill $\Box$

\vspace{.2cm}
\noindent {\bf Note.} This manuscript dates back to the years 2010--2011, partly as an offspring of the work \cite{gowda-sznajder}, and contains essentially one result. It is now published with hope that it would support some bigger piece of work. Due to the increased interest in applications of the Lorentz cone, the author has decided to record this in arXiv.org. 

\vspace{.2cm}   
\noindent {\bf Acknowledgements.} \,The author is indebted to Professor M. Seetharama Gowda for suggestions that led to a final version of this paper. 


\end{document}